\documentclass[letterpaper, 10 pt, conference]{ieeeconf} 
\IEEEoverridecommandlockouts                              % This command is only
\usepackage{graphics} % for pdf, bitmapped graphics files
\usepackage{epsfig} % for postscript graphics files
\usepackage{cite}

\usepackage{amsmath} % assumes amsmath package installed
\usepackage{amssymb}  % assumes amsmath package installed
\usepackage{subcaption}

\usepackage{mathrsfs}

\title{\LARGE \bf
A Closed-Form Analytical Solution for Optimal Coordination of Connected and Automated Vehicles}

\author{Andreas A. Malikopoulos, {\itshape{Senior Member, IEEE}}, and Liuhui Zhao, {\itshape{Member, IEEE}}
\thanks{This research was supported in part by the  Delaware Energy Institute (DEI) and by ARPAE's NEXTCAR program under the award number DE-AR0000796.}%
\thanks{The authors are with the Department of Mechanical Engineering, University of Delaware, Newark, DE 19716 USA (email: andreas@udel.edu; lhzhao@udel.edu.} }

\begin{document}

\maketitle
\thispagestyle{empty}
\pagestyle{empty}

%%%%%%%%%%%%%%%%%%%%%%%%%%%%%%%%%%%%%%%%%%%%%%%%%%%%%%%%%%%%%%%%%%%%%%%%%%%%%%%%
\begin{abstract}
In earlier work, a decentralized optimal control framework was established for coordinating online connected and automated vehicles (CAVs) in merging roadways, urban intersections, speed reduction zones, and roundabouts. The dynamics of each vehicle were represented by a double integrator and the Hamiltonian analysis was applied to derive an analytical solution that minimizes the $L^2$-norm of the control input. However, the analytical solution did not consider the rear-end collision avoidance constraint. In this paper, we derive a complete, closed-form analytical solution that includes the rear-end safety constraint in addition to the state and control constraints. We augment the double integrator model that represents a vehicle with an additional state corresponding to the distance from its preceding vehicle. Thus, the rear-end collision avoidance constraint is included as a state constraint. The effectiveness of the solution is illustrated through simulation.
\end{abstract}

%%%%%%%%%%%%%%%%%%%%%%%%%%%%%%%%%%%%%%%%%%%%%%%%%%%%%%%%%%%%%%%%%%%%%%%%%%%%%%%%
\indent

%\begin{IEEEkeywords}
%Connected and automated vehicles, decentralized optimal control, autonomous intersections, traffic flow,  energy usage, safety.
%\end{IEEEkeywords}
%%%%%%%%%%%%%%%%%%%%%%%%%%%%%%%%%%%%
%SECTION I: INTRODUCTION
%%%%%%%%%%%%%%%%%%%%%%%%%%%%%%%%%%%%

%\pagestyle{plain} ----------- FOR PAGES IN THE PAPER %%%%%%%%%%%%%%%%%%%%%%%%%%%%%%%%%%%%%%%%%%%%%%

\section{Introduction} \label{sec:1}
We are currently witnessing an increasing integration of our energy, transportation, and cyber networks, which, coupled with the human interactions, is giving rise to a new level of complexity in the transportation network. As we move to increasingly complex \cite{Malikopoulos2016b} emerging mobility systems, new control approaches are needed to optimize the impact on system behavior of the interplay between vehicles at different transportation scenarios, e.g., intersections, merging roadways, roundabouts, speed reduction zones. These scenarios along with the driver responses to various disturbances are the primary sources of bottlenecks that contribute to traffic congestion.
More recently, a study \cite{Ratti2016} indicated that transitioning from intersections with traffic lights to autonomous intersections, where vehicles can coordinate and cross the intersection without the use of traffic lights, has the potential of doubling capacity and reducing delays. 

Several research efforts have been reported in the literature proposing either centralized or decentralized approaches on coordinating CAVs at intersections. Dresner and Stone \cite{Dresner2004} proposed the use of the reservation scheme to control a single intersection of two roads with vehicles traveling with similar speed on a single direction on each road. Some approaches have focused on coordinating vehicles at intersections to improve the travel time. Kim and Kumar \cite{Kim2014} proposed an approach based on model predictive control that allows each vehicle to optimize its movement locally in a distributed manner with respect to any objective of interest. Colombo and Del Vecchio \cite{Colombo2014} constructed the invariant set for the control inputs that ensure lateral collision avoidance. Previous work has also focused on multi-objective optimization problems for intersection coordination, mostly solved as a receding horizon control problem, in either centralized or decentralized approaches \cite{Kamal2013a, Kamal2014, Campos2013, Kim2014, qian2015}. 
For instance, Campos et al. \cite{Campos2014} applied a receding horizon framework for a decentralized solution for autonomous vehicles driving through traffic intersections. Qian et al. \cite{qian2015} proposed to solve the intersection coordination problem in two levels, where vehicles coordination was handled based on predefined priority scheme at the upper level, and each vehicle solved its own multi-objective optimization problem at the lower level.  
A detailed discussion of the research efforts in this area that have been reported in the literature to date can be found in \cite{Malikopoulos2016a}. 

Coordinating CAVs at an urban intersection generally involves a two-level joint optimization problem: (1) an upper level vehicle coordination problem which specifies the sequence that each CAV crosses the intersection \cite{Malikopoulos2019a} and (2) a lower level optimal control problem in which each CAV derives its optimal acceleration/deceleration, in terms of energy, to cross the intersection. In earlier work, a decentralized optimal control framework was established for coordinating online CAVs in different transportation scenarios, e.g., merging roadways, urban intersections, speed reduction zones, and roundabouts. The analytical solution using a double integrator model, without considering state and control constraints, was presented in \cite{Rios-Torres2015}, \cite{Rios-Torres2}, and \cite{Ntousakis:2016aa} for coordinating online CAVs at highway on-ramps, in \cite{Zhang2016a} at two adjacent intersections, and in \cite{Malikopoulos2018a} at roundabouts. The solution of the unconstrained problem was also validated experimentally at the University of Delaware's Scaled Smart City using 10 robotic cars \cite{Malikopoulos2018b} in a merging roadway scenario. The solution of the optimal control problem considering state and control constraints was presented in \cite{Malikopoulos2017} at an urban intersection, without considering rear-end collision avoidance constraint though. The conditions under which the rear-end collision avoidance constraint never becomes active were discussed in \cite{Malikopoulos2018c}.

In this paper, we consider that the sequence that each CAV crosses the intersection is given and we focus only on the lower level optimal control problem. We derive a complete, closed-form analytical solution that includes the rear-end safety constraint in addition to the state and control constraints of the lower level problem. We augment the double integrator model that represents a vehicle with an additional state corresponding to the distance from its preceding vehicle. Thus, the rear-end collision avoidance constraint is included as a state constraint. Furthermore, we allow the safe distance between two vehicles to be a function of the vehicle's speed. 

The structure of the paper is organized as follows. In Section II, we review the problem of vehicle coordination at an urban intersection and provide the modeling framework. In Section III, we derive the analytical, closed form solution. In Section IV, we validate the effectiveness of the analytical solution through simple driving scenarios. Finally, we offer concluding remarks in Section V.

%%%%%%%%%%%%%%%%%%%%%%%%%%%%%%%%%%%%
%%%%%%%%%%%%%%%%%%%%%%%%%%%%%%%%%%%%
%SECTION II: Problem Formulation
%%%%%%%%%%%%%%%%%%%%%%%%%%%%%%%%%%%%
%%%%%%%%%%%%%%%%%%%%%%%%%%%%%%%%%%%%
\section{Problem Formulation} \label{sec:2}
\subsection{Vehicle Model, Constraints, and Assumptions} \label{sec:2a}
\begin{figure}[ptb]
\centering
\includegraphics[width=2.4 in]{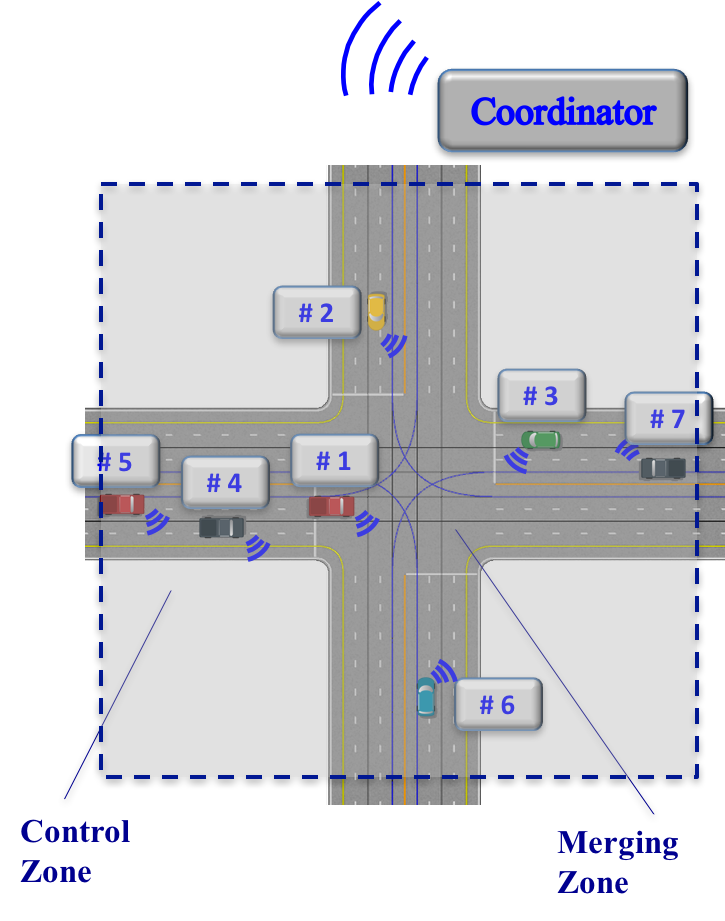} \caption{An urban intersection with
connected and automated vehicles.}%
\label{fig:1}%
\end{figure}

We consider a single urban intersection (Fig. \ref{fig:1}). The region at the center of the intersection, called \textit{merging zone}, is the area of potential lateral collision of the vehicles. The intersection has a \textit{control zone} and a coordinator that can communicate with the vehicles traveling inside the control zone. Note that the coordinator is not involved in any decision on the vehicle. The distance from the entry of  the control zone until the entry of the merging zone is $L$, and it is assumed to be the same for all entry points of the control zone. Note that the $L$ could be in the order of hundreds of $m$ depending on the coordinator's communication range capability, while $S$ is the length of a typical intersection.

Let $N(t)\in\mathbb{N}$ be the number of CAVs inside the control zone at time
$t\in\mathbb{R}^{+}$ and $\mathcal{N}(t)=\{1,\ldots,N(t)\}$ be a queue which
designates the order in which these vehicles will be entering the merging zone. Let $t_{i}^{f}$ be the assigned time for vehicle $i$ to exits the control zone.
There is a number of ways to assign $t_{i}^{f}$ for each CAV $i$. For example, we may
impose a strict first-in-first-out queueing structure, where each vehicle must
enter the merging zone in the same order it entered the control zone. The policy through which the \textquotedblleft schedule\textquotedblright~  is specified is the result of a higher level optimization problem. This policy, which  determines the time $t_{i}^{f}$ that each CAV $i$ exits the control zone, can aim at maximizing the throughput at the intersection while ensuring that the lateral collision avoidance constraint never becomes active. Once the desired $t_{i}^{f}$ for each CAV $i$ is determined, it is stored in the coordinator and is not changed.  On the other hand, for each CAV $i$, deriving the optimal control input (minimum acceleration/deceleration) to achieve the target $t_{i}^{f}$ can aim at minimizing its fuel consumption \cite{Malikopoulos2010a} while ensuring that the rear-end collision avoidance constraint never becomes active.

In what follows, we assume that a scheme for determining
$t_{i}^{f}$ (upon arrival of CAV $i$) is given, and we will focus on a lower level control problem that will yield for each CAV the optimal control input (acceleration/deceleration) to achieve the assigned $t_{i}^{m}$ subject to the state, control, and rear-end collision avoidance constraints.

\subsection{Vehicle Model and Constraints}

We consider a number of CAVs $N(t)\in\mathbb{N}$, where $t\in\mathbb{R}$ is the time, that enter the control zone. 
%Let $\mathcal{N}(t)=\{1,\ldots,N(t)\},$ be a queue of the CAVs
%associated with the control zone of the intersection. Each vehicle has a unique location $i$ in the  queue, where $i=N(t)+1$. 
We represent the dynamics of each vehicle
$i\in\mathcal{N}(t)$, 
with a state equation
\begin{equation}
\dot{x_{i}}=f(t,x_{i},u_{i}),\qquad x_{i}(t_{i}^{0})=x_{i}^{0},\label{eq:model}\\
\end{equation}
where $t\in\mathbb{R}^{+}$ is the time, $x_{i}(t)$, $u_{i}(t)$ are the state
of the vehicle and control input, $t_{i}^{0}$ is the time that vehicle $i$
enters the control zone, and $x_{i}^{0}$ is the value of the initial
state. We assume that the dynamics of each vehicle are
\begin{equation}%
\begin{split}
\dot{p}_{i} &  =v_{i}(t)\\
\dot{v}_{i} &  =u_{i}(t)\\
\dot{s}_{i} &  = \xi_i \cdot (v_{k}(t)-v_{i}(t))
\label{eq:model2}
\end{split}
\end{equation}
where $p_{i}(t)\in\mathcal{P}_{i}$, $v_{i}(t)\in\mathcal{V}_{i}$, and
$u_{i}(t)\in\mathcal{U}_{i}$ denote the position, speed and
acceleration/deceleration (control input) of each vehicle $i$ inside the control zone; $s_{i}(t)\in\mathcal{S}_{i}, s_{i}(t)=p_{k}(t)-p_{i}(t) $ denotes the distance of vehicle $i$ from the vehicle $k$ which is physically immediately ahead of $i$, and $\xi_{i}$ is a reaction constant of the vehicle. The sets $\mathcal{P}_{i}$,
$\mathcal{V}_{i}$, $\mathcal{U}_{i}$, and $\mathcal{S}_{i}$, $i\in\mathcal{N}(t),$
are complete and totally bounded subsets of $\mathbb{R}$.

Let
$x_{i}(t)=\left[p_{i}(t) ~ v_{i}(t) ~ s_{i}(t)\right]  ^{T}$ denote the state of each vehicle $i$, with initial value
$x_{i}^{0}=\left[p_{i}^{0} ~ v_{i}^{0} ~s_{i}^{0}\right]  ^{T}$, where $p_{i}^{0}= p_{i}(t_{i}^{0})=0$ at the entry of the control zone, taking values in $\mathcal{X}_{i}%
=\mathcal{P}_{i}\times\mathcal{V}_{i}$.  The state space 
$\mathcal{X}_{i}$ for each vehicle $i$ is
closed with respect to the induced topology on $\mathcal{P}_{i}\times
\mathcal{V}_{i}$ and thus, it is compact.
We need to ensure that for any initial state $(t_i^0, x_i^0)$ and every admissible control $u(t)$, the system \eqref{eq:model} has a unique solution $x(t)$ on some interval $[t_i^0, t_i^f]$, where $t_i^f$ is the time that vehicle $i\in\mathcal{N}(t)$ exits the \textit{control zone}. 
The following observations from \eqref{eq:model} satisfy some regularity conditions required both on $f$ and admissible controls $u(t)$ to guarantee local existence and uniqueness of solutions for \eqref{eq:model}: a) The function $f$ is continuous in $u$ and continuously differentiable in the state $x$, b) The first derivative of $f$ in $x$, $f_x$, is continuous in $u$, and c) The admissible control $u(t)$ is continuous with respect to $t$.

To ensure that the control input and vehicle speed are within a
given admissible range, the following constraints are imposed.
\begin{equation}%
\begin{split}
u_{i,min} &  \leq u_{i}(t)\leq u_{i,max},\quad\text{and}\\
0 &  \leq v_{min}\leq v_{i}(t)\leq v_{max},\quad\forall t\in\lbrack t_{i}%
^{0},t_{i}^{f}],
\end{split}
\label{speed_accel constraints}%
\end{equation}
where $u_{i,min}$, $u_{i,max}$ are the minimum deceleration and maximum
acceleration for each vehicle $i\in\mathcal{N}(t)$, and $v_{min}$, $v_{max}$ are the minimum and maximum speed limits respectively. 

To ensure the absence of rear-end collision of two consecutive vehicles traveling on the same lane,  the position of the preceding vehicle should be greater than or equal to the position of the following vehicle plus a predefined safe distance $\delta_i(t)$. Thus we impose the rear-end safety constraint 
\begin{equation}
\begin{split}
s_{i}(t)=\xi_i \cdot (p_{k}(t)-p_{i}(t)) \ge \delta_i(t),~ \forall t\in [t_i^0, t_i^f],
\label{eq:rearend}
\end{split}
\end{equation}
where $k$ is some vehicle which is physically immediately ahead of $i$ in the same lane. We relate the minimum safe distance $\delta_i(t)$ as a function of speed $v_i(t)$, 
\begin{equation}
\begin{split}
\delta_i(t)=\gamma_i + \rho_i \cdot v_i(t),~ \forall t\in [t_i^0, t_i^f],
\label{eq:safedist}
\end{split}
\end{equation}
where $\gamma_i$ is the standstill distance, and $\rho_i$ is minimum time gap that vehicle $i$ would maintain while following another vehicle.

Once the time $t_{i}^{f}$ that each vehicle $i\in\mathcal{N}(t)$ will be exiting the control zone is assigned, the problem for each vehicle is to minimize the cost functional $J_{i}(u(t))$, which is the $L^2$-norm of the control input in $[t_i^0, t_i^f]$
\begin{gather}\label{eq:decentral}
\min_{u(t)\in U_i} J_{i}(u(t))=  \frac{1}{2} \int_{t^0_i}^{t^f_i} u^2_i(t)~dt,\\ 
\text{subject to}%
:\eqref{eq:model2},\eqref{speed_accel constraints},\eqref{eq:rearend},\nonumber\\
\text{and given }t_{i}^{0}\text{, }v_{i}^{0}\text{, } s_i(t_i^0) = s_i^0\text{, } p_{i}(t_{i}^{0})=0\text{, } p_{i}(t_{i}^{f})=L+S, t_{i}^{f}.\nonumber
\end{gather}

%%%%%%%%%%%%%%%%%%%%%%%%%%%%%%%%%%%%
%%%%%%%%%%%%%%%%%%%%%%%%%%%%%%%%%%%%
%SECTION III: Analytical solution of the decentralized control problem
%%%%%%%%%%%%%%%%%%%%%%%%%%%%%%%%%%%%
%%%%%%%%%%%%%%%%%%%%%%%%%%%%%%%%%%%%
\section{Analytical solution of the optimal control problem} \label{sec:3}

Let $\mathbf{S}_i(t,x(t))$ be the vector of the constraints in \eqref{eq:decentral} which do not explicitly depend on $u(t)$ \cite{bryson1975applied},
\begin{gather} 
\mathbf{S}_{i}\big(t,x(t)\big) = % \left\{
\begin{bmatrix}
v_{i}(t) - v_{max} \\
v_{min} - v_{i}(t) \\
\delta_{i}(t) - s_{i}(t)
\end{bmatrix}. \notag 
 \label{eq:gNoU}
\end{gather}
Since $\mathbf{S}_i(t)\leq \mathbf{0}$ is satisfied for all $t\in [t_i^0, t_i^f]$, it follows that $\dot{\mathbf{S}}(t) =  \leq \mathbf{0}$. 

Thus, the  Hamiltonian becomes
\begin{gather}
H_{i}\big(t, p_{i}(t), v_{i}(t), s_{i}(t), u_{i}(t)\big)  \nonumber \\
=\frac{1}{2} u(t)^{2}_{i} + \lambda^{p}_{i} \cdot v_{i}(t) + \lambda^{v}_{i} \cdot u_{i}(t) +\lambda^{s}_{i} \cdot \xi_i \cdot (v_{k}(t) - v_{i}(t)) \nonumber\\
+ \mu^{a}_{i} \cdot(u_{i}(t) - u_{max})
+ \mu^{b}_{i} \cdot(u_{min} - u_{i}(t)) \nonumber\\
+ \mu^{c}_{i} \cdot  u_{i}(t) - \mu^{d}_{i} \cdot u_{i}(t) \nonumber\\ 
+ \mu^{s}_{i} \cdot (\rho_i \cdot u_i(t) - \xi_i\big(v_{k}(t) - v_i(t)\big)) ,\label{eq:16b}
\end{gather}
where $\lambda^{p}_{i}$, $\lambda^{v}_{i}$, and $\lambda^{s}_{i}$ are the influence functions \cite{bryson1975applied}, and
$\mu^{T}$ is the vector of the Lagrange multipliers.

For each $i\in\mathcal{N}(t)$, the Euler-Lagrange equations are
\begin{gather}\label{eq:EL1}
\dot\lambda^{p}_{i}(t) = - \frac{\partial H_i}{\partial p_{i}} = 0, \\
\dot\lambda^{v}_{i}(t) = - \frac{\partial H_i}{\partial v_{i}} = - (\lambda^{p}_{i} - \lambda^{s}_{i} \cdot \xi_i + \mu^s_i \cdot \xi_i),  \label{eq:EL2}\\
\dot\lambda^{s}_{i}(t) = - \frac{\partial H_i}{\partial s_{i}} = 0, \label{eq:EL3}%
\end{gather}
\begin{equation}
\label{eq:KKT1}
\frac{\partial H_i}{\partial u_{i}} = u_{i}(t) + \lambda
^{v}_{i} + \mu^{a}_{i} - \mu^{b}_{i} + \mu^c_i - \mu^d_i + \mu^s_i \rho_i = 0,
\end{equation}
with boundary conditions 
\begin{gather}
p_i(t_i^0) = p_i^0, ~ v_i(t_i^0) = v_i^0, ~ s_i(t_i^0) = s_i^0,\notag \\
p_i(t_i^f) = L+S, ~ \lambda^v_i(t_i^f) = 0, ~ \lambda^s_i(t_i^f) = 0,
\label{eq:bound}
\end{gather}
where $\lambda^v_i(t_i^f)=\lambda^s_i(t_i^f) = 0$ since the states $v_i(t_i^f)$ and $s_i(t_i^f)$ are not prescribed at $t_i^f$ \cite{bryson1975applied}.

To address this problem, the constrained and unconstrained arcs will be pieced together to satisfy the Euler-Lagrange equations and necessary condition of optimality. Based on our state and control constraints \eqref{speed_accel constraints}, \eqref{eq:rearend} and boundary conditions, the optimal solution is the result of different combinations of the following possible arcs.

%%%%%%%%%%%%%%%%%%%%%%%%%%%%%%%%%%
%Case 1
%%%%%%%%%%%%%%%%%%%%%%%%%%%%%%%%%%

\subsubsection{Inequality State and Control Constraints are not Active}

In this case, we have 
$\mu^{a}_{i} = \mu^{b}_{i}= \mu^{c}_{i}=\mu^{d}_{i}=\mu^{e}_{i}=0.$
Applying the necessary condition \eqref{eq:KKT1}, the
optimal control can be given 
\begin{equation}
u_{i}(t) + \lambda^{v}_{i}= 0, \quad i \in\mathcal{N}(t). \label{eq:17}
\end{equation}
From (\ref{eq:EL1}), \eqref{eq:EL2}, and \eqref{eq:EL3}  we have $\lambda^{p}_{i}(t) = a_{i}$,   
$\lambda^{s}_{i}(t)= b_{i}$, and $\lambda^{v}_{i}(t) = -\big((a_{i} - b_{i} \cdot \xi_i)\cdot t + c_{i}\big)$. 
The coefficients $a_{i}$, $b_{i}$, and $c_{i}$
are constants of integration corresponding to each vehicle $i$. From \eqref{eq:17}
the optimal control input (acceleration/deceleration) as a function of time is given by
\begin{equation}
u^{*}_{i}(t) = (a_{i} - b_{i} \cdot \xi_i) \cdot t + c_{i}, ~ \forall t \ge t^{0}_{i}. \label{eq:20}
\end{equation}

Substituting the last equation into \eqref{eq:model2} we find the optimal speed and position for each vehicle,
namely
\begin{gather}
v^{*}_{i}(t) = \frac{1}{2} (a_{i} - b_{i} \cdot \xi_i) \cdot t^2 + c_{i} \cdot t +d_{i}, ~ \forall t \ge t^{0}_{i}, \label{eq:21}\\
p^{*}_{i}(t) = \frac{1}{6} (a_{i} - b_{i} \cdot \xi_i) \cdot t^3 +\frac{1}{2} c_{i} \cdot t^2 + d_{i}\cdot t +e_{i}, ~ \forall t \ge t^{0}_{i}, \label{eq:22}%
\end{gather}
where $d_{i}$ and $e_{i}$ are constants of integration. The constants of integration $a_i$, $c_{i}$, $d_{i}$, and $e_{i}$ are computed at each time $t, t^{0}_{i} \le t \le t^{f}_{i}$, using the values of the control input, speed, and position of each vehicle $i$ at $t$, the position $p_{i}(t^{f}_{i})$, and the values of the one of terminal transversality condition, i.e., $\lambda^{v}_{i}(t^{f}_{i})$. Since the terminal cost, i.e., the control input, at $t^{f}_{i}$ is zero,  we can assign $\lambda^{v}_{i}(t^{f}_{i}) =0$. 
%To
%derive online the optimal control for each vehicle $i$, we need to update the integration constants at each time $t$. We form the following system of six equations, namely%
%\begin{gather}
%\left[
%\begin{array}
%[c]{ccccc}%
%\frac{t^2}{2} & -\frac{t^2}{2} \cdot \xi_i & t & 1 & 0 \\
%\frac{t^3}{6} & -\frac{t^3}{6} \cdot \xi_i & \frac{t^2}{2} & t & 1\\
%\frac{(t^{f}_{i})^3}{6} & -\frac{(t^{f}_{i})^3}{6} \cdot \xi_i & \frac{(t^{f}_{i})^2}{2} & t^{f}_{i} & 1\\
%-t_i^f & t_i^f\cdot \xi_i & -1 & 0 & 0 \\
%0 & 1 & 0 & 0 & 0\\
%\end{array}
%\right] \cdot
%\left[
%\begin{array}
%[c]{c}%
%a_{i}\\
%b_{i}\\
%c_{i}\\
%d_{i}\\
%e_{i}
%\end{array}
%\right] \nonumber \\ 
%=\left[
%\begin{array}
%[c]{c}%
%v_{i}(t)\\
%p_{i}(t)\\
%p_{i}(t^{f}_{i})\\
%\lambda^{v}_{i}(t^{f}_{i})\\
%\lambda^{s}_{i}(t^{f}_{i})
%\end{array}
%\right], \forall t \ge t^{0}_{i}.\label{eq:23}%
%\end{gather}
%Since \eqref{eq:23} is computed online, the controller yields the
%optimal control online for each vehicle $i$, with feedback provided through the re-calculation of the constants of integration 
%$\tilde{a_i}$, $c_{i}$, $d_{i}$, and $e_{i}$.

%%%%%%%%%%%%%%%%%%%%%%%%%%%%%%%%%%
%Case 1a, s_i(t) = \detla
%%%%%%%%%%%%%%%%%%%%%%%%%%%%%%%%%%

\subsubsection{The State Constraint $s_{i}(t)=\delta(t)$ Becomes Active}
\label{sec:3b}
Suppose vehicle $i\in\mathcal{N}(t)$ starts from a feasible state and control at $t=t_{i}^0$ and at some time $t=t_{1}\le t_i^f$, $s_{i}(t_{1})=\delta(t_{1})$ while $v_{min} < v_{i}(t_{1}) < v_{max}$ and $u_{i,min} < u_{i}(t_{1}) < u_{i,max}$. In this case, $\mu_{i}^{s} \neq 0$. 

Let $N_i(t,x(t))= \gamma_i+ \rho_i v_i^*(t_1) - \xi_ip_k^*(t_1)+\xi_ip_i^*(t_1).$ Then, we have
\begin{equation} \label{eq:delta1}
N_i(t,x(t))=  \gamma_i+ \rho_i v_i^*(t_1) - \xi_i (p_k^*(t_1)+ p_i^*(t_1)) = 0,
\end{equation}
which represents a terminal constraint for the state $s_i(t)$ in $t\in [t_{i}^0, t_1].$ 
Since $N_i(t_1,x(t_1))=0$, its first derivative should vanish
\begin{equation} \label{eq:delta2}
\dot{N_i}(t_1,x(t_1))=  \rho_i u_i^*(t_1) - \xi_i (v_k^*(t_1) - v_i^*(t_1)) = 0,
\end{equation}
from which we derive the value of the optimal control at $t=t_1^+$
\begin{equation} \label{eq:delta3}
u_i^*(t_1^+) = \frac{ \xi_i (v_k^*(t_1^+) - v_i^*(t_1^+))}{\rho_i}.
\end{equation}

From \eqref{eq:20} and \eqref{eq:delta3}, we note that the optimal control input is not continuous at $t_1$, hence the junction point at $t_1$ is a corner.
The boundary conditions at the corner for the influence fundtions are
\begin{equation}
\lambda_i^T(t_1^-)=\lambda_i^T(t_1^+)+\pi_i \frac{\partial{N_i(t,x(t))}}{\partial x(t)}.
\label{eq:delta4}
\end{equation}
The transversality condition is
\begin{equation}
\lambda_i^T(t_1^-) \dot{x}(t_1^-)=\lambda_i^T(t_1^+) \dot{x}(t_1^+) - \pi_i \frac{\partial{N_i(t,x(t))}}{\partial t_1},
\label{eq:trans}
\end{equation}
where $\pi_i$ is a  Lagrange multiplier constant. The influence functions, $\lambda_i^T(t_1^+)$, at $t_1^+$, the entry time  $t_1$, and the Lagrange multipllier $\pi_i$ constitute 3+1+1 quantities that are determined so as to satisfy \eqref{eq:delta1}, \eqref{eq:delta4}, and \eqref{eq:trans}.  Note, the values of the influence functions $\lambda_i^T(t_1^-)$ at $t_1^-$ are known from the unconstrained arc in $[t_i^0,t_1]$, i.e., $\lambda_i^p(t_1^-)=\alpha_i$, $\lambda_i^v(t_1^-)= -u_i^*(t_1^-),$ $\lambda_i^s(t_1^-)=\beta_i$, and the state variables are continuous at the junction point, $t_1$, i.e., $p_i(t_1^-) = p_i(t_1^+),$ $v_i(t_1^-) = v_i(t_1^+),$ $s_i(t_1^-) = s_i(t_1^+).$ The unconstrained and constrained arcs are pieced together to determine the 3+1+1 quantities above along with the constants of integration in \eqref{eq:20}-\eqref{eq:22} while the Hamiltonian at the corner is $H_i(t_1^-)=H_i(t_1^+)- \pi_i \frac{\partial{N_i(t,x(t))}}{\partial t_1}$.

For the optimal control of the constrained arc, $\delta_i(t)-s_i(t) \le 0$, we have the following two cases to consider: (a) when the speed, $v_k(t),$ of the preceding vehicle $k$ is decreasing, and (b) when the speed, $v_k(t),$ of the preceding vehicle $k$ is either increasing or it is constant.

\paragraph{The speed, $v_k(t),$ of the preceding vehicle $k$ is decreasing}

\textit{Case 1: The exit point $t_2$ leads to the arc $u_{i,\min}-u_i(t)\le 0$.} Next, we consider the case that the exit point $t_2$ of the constrained arc, $\delta_i(t)-s_i(t) \le 0,$ leads to the arc $u_{i,\min}-u_i(t)\le 0$. It follows that
\begin{equation}
u_i^*(t) = u_{i,\min}, ~t\in[t_2,t_3], 
\label{eq:56}
\end{equation}
where $t_3$ is the exit point of the arc $u_{i,\min}-u_i(t)\le 0$.
By integrating \eqref{eq:56} we have
\begin{align}
v_i^*(t) &= u_{i,\min}\cdot t + h_i, ~t\in[t_2,t_3],\\
p_i^*(t) &= u_{i,\min} \cdot \frac{t^2}{2} + h_i t + q_i, ~t\in[t_2,t_3],
\end{align}
where $h_i$ and $q_i$ are constants of integration. For the exit point $t_3$ of this arc, we have the following result.

If vehicle $i$ remains at the constrained arc, $u_{i,\min}-u_i \le 0$, until $t_i^f$, then we use the interior constraints and boundary conditions from which we can compute $t_2$, and the constants of integrations $h_i$ and $q_i$. If, however, at some time $t=t_3$, vehicle $i$ exits the constrained arc, $u_{i,\min}-u_i(t) \le 0$, and enters the arc $v_{\min} - v_i(t) \le 0$, then it follows that $u_{i}^*(t) = 0,$ for all $t\in[t_3, t_i^f]$, and the optimal speed and position of $i$ are
\begin{align}
v^{*}_{i}(t) &=  v_{min}\label{eq:cas5a2}, ~ t\in[t_3, t_i^f], \\
p^{*}_{i}(t) &= v_{min}~t + r_{i} \label{eq:cas5a3}, ~ t\in[t_3, t_i^f],
\end{align}
where $r_{i}$ is a constant of integration. In this case, we piece together the unconstrained with the two constrained arcs, $u_{i,\min}-u_i(t) \le 0$ and $v_{\min} - v_i(t) \le 0$, to satisfy the interior constraints and boundary conditions from which we can compute $t_2$, $t_3$ and the constants of integration $h_i, q_i,$ and $r_i$.

\textit{Case 2: The exit point $t_2$ leads to the arc $v_{\min}-v_i(t)\le 0$.} Next, we consider the case that the exit point $t_2$ of the constrained arc, $\delta_i(t)-s_i(t) \le 0,$ leads to the arc $v_{\min}-v_i(t)\le 0$. It follows that $u_{i}^*(t) = 0,$ for all $t\in[t_2, t_i^f]$, and the optimal speed and position of the vehicle are given by \eqref{eq:cas5a2} and \eqref{eq:cas5a3}. From the interior constraints and boundary conditions, we can compute $t_2$, and the constant of integration $r_i$.

\paragraph{The speed, $v_k(t),$ of the preceding vehicle $k$ is either increasing or it is constant}

The unconstrained arc for all  $t\in[t_2, t_i^f],$ consists of a set of equations as in \eqref{eq:20} - \eqref{eq:22} for the optimal control, speed, and position of vehicle $i$, i.e.,  $u^{*}_{i}(t) = a_i'\cdot t + c_{i}',$ $v^{*}_{i}(t) = \frac{1}{2} a_i' \cdot t^2 + c_{i}' \cdot t +d_{i}',$ and $p^{*}_{i}(t) = \frac{1}{6} a_i' \cdot t^3 +\frac{1}{2} c_{i}' \cdot t^2 + d_{i}'\cdot t +e_{i}',$ where $a_i',$ $c_i',$ $d_i',$ and $e_i',$ are constants of integration that can be computed along with $t_2$ from the interior constraints and boundary conditions.

Similar results are obtained for the remaining cases. To derive the analytical solution of \eqref{eq:decentral}, we first start with the unconstrained arc and derive the solution using  \eqref{eq:20} - \eqref{eq:22}. If the solution violates any of the state or control constraints, then the unconstrained arc is pieced together with the arc corresponding to the violated constraint, and we re-solve the problem with the two arcs pieced together. The two arcs yield a set of algebraic equations which are solved simultaneously using the boundary conditions of \eqref{eq:decentral} and interior conditions between the arcs. If the resulting solution, which includes the determination of the optimal switching time from one arc to the next one, violates another constraint, then the last two arcs are pieced together with the arc corresponding to the new violated constraint, and we re-solve the problem with the three arcs pieced together. The three arcs will yield a new set of algebraic equations that need to be solved simultaneously using the boundary conditions of \eqref{eq:decentral} and interior conditions between the arcs. The resulting solution includes the optimal switching time from one arc to the next one. The process is repeated until the solution does not violate any other constraints.

%%%%%%%%%%%%%%%%%%%%%%%%%%%%%%%%%%%%%%%%%%%%%%%%%%%%%%%%%%%%%%%%%%%%%%%%%%%%%%%%
%%%%%%%%%%%%%%%%%%%%%%%%%%%%%%%%%%%%
%SECTION IV: Simulation Results
%%%%%%%%%%%%%%%%%%%%%%%%%%%%%%%%%%%%
\section{Simulation Results}
To validate the effectiveness of the analytical solution for real-end collision avoidance, we created a simple driving scenario in MATLAB. The length of the control zone is 300 $m$. The following vehicle $i$ is located at the entry of the control zone (Fig. \ref{fig:1}) with the initial speed of 14 $m/s$. At the time that vehicle $i$ enters the control zone, the leading vehicle $k$ has a speed of 11.5 $m/s$ and is located at 20 $m$ (inside the control zone). In this analysis, we set -1 $m/s^2$ and 1 $m/s^2$ as the minimum and maximum acceleration. For simplification, we set the final time for vehicle $i$ is 26 $s$. We analyzed three cases with different leading vehicle acceleration profiles to test the effectiveness of our model. For comparison, we also include the scenario when the safety constraint is not considered in the optimization model.

\subsection{Case 1: constant acceleration of leading vehicle}
\begin{figure}[ptb]
\centering
\begin{subfigure}[b]{0.45\textwidth}
\centering
\includegraphics[width=\textwidth]{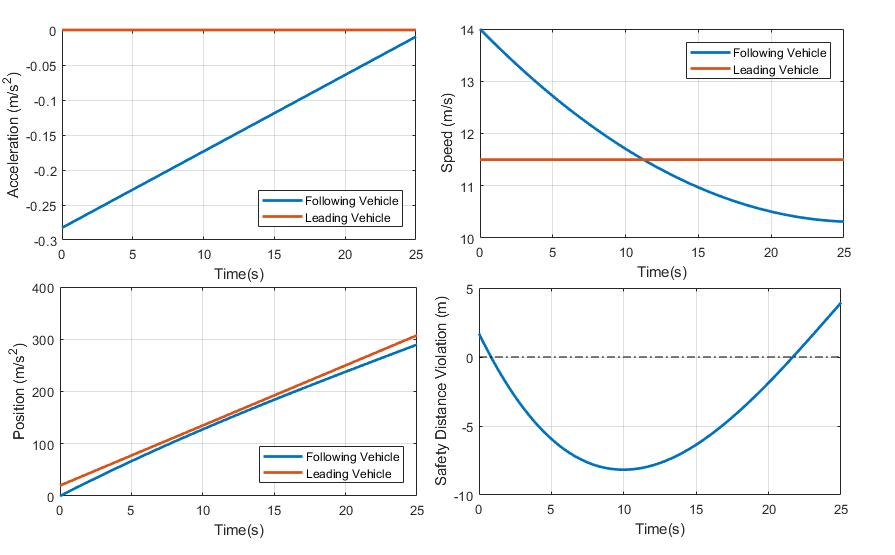} 
\caption{Results without safety constraint} \label{fig:1a}
\end{subfigure}

\begin{subfigure}[b]{0.45\textwidth}
\centering
\includegraphics[width=\textwidth]{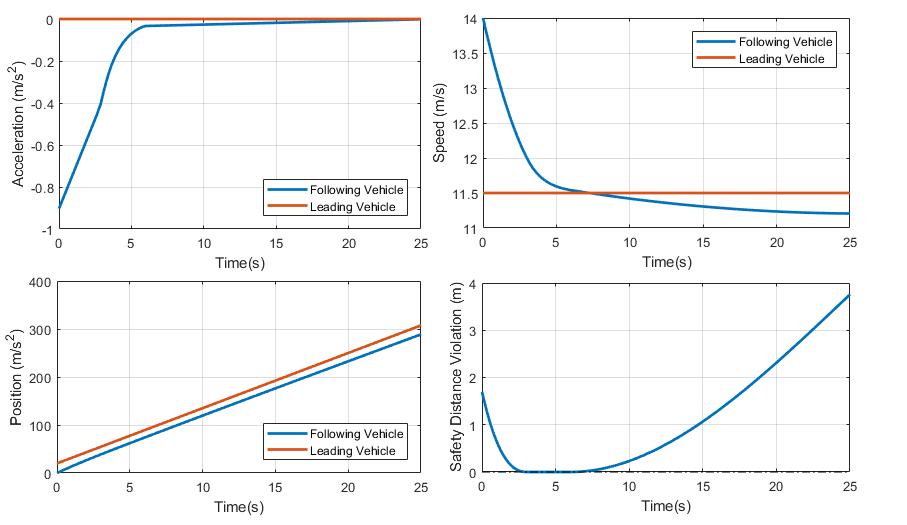} 
\caption{Results with safety constraint} \label{fig:1b}
\end{subfigure}

\caption{Optimization results for case 1.}
\label{fig:case1}%
\end{figure}

In this case, we consider constant speed of leading vehicle $k$. We see that in Fig. \ref{fig:1a}, if the safety constraint is not incorporated in the optimization model, linear acceleration profile is yielded, however, the following distance of vehicle $i$ violates the minimum safety distance. Two vehicles get too close to each other, which creates an extremely unsafe driving situation. Considering the safety constraint, the optimal acceleration profile is presented in Fig. \ref{fig:1b}. We observe three arcs in the optimal acceleration profile, before $t_1=3.1$ $s$, the safety constraint is not violated, vehicle $i$ decelerates with a much lower acceleration than the recommended acceleration without safety constraint. At $t_1=3.1$ $s$, safety constraint is violated, vehicle $i$ enters the constrained arc at $t_1=3.1$ $s$ and leaves constrained arc at $t_2=6.5$ $s$.

\subsection{Case 2: linearly decreasing acceleration of leading vehicle}
\begin{figure}[ptb]
	\centering
	\includegraphics[width=0.45\textwidth]{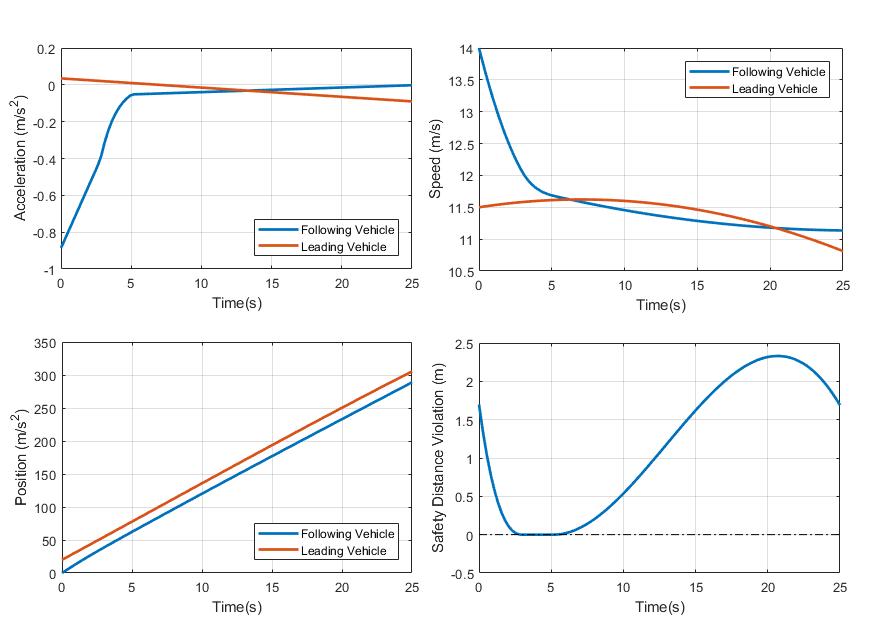}   
	\caption{Optimization results for case 2.} 
	\label{fig:case2}
\end{figure}

In this case, we consider a decreasing acceleration profile of vehicle $k$ with a positive initial acceleration. Similar to case 1, the safety constraint is activated. By piecing together the unconstrained and constrained arcs, the results corresponding to the closed form analytical solution are shown in Fig. \ref{fig:case2}. Before the entry time at $t_1=2.9$ $s$, vehicle $i$ travels with a linearly decreasing acceleration until the safety constraint is activated (i.e., $s_i(t)-\delta(t)=0$). Since vehicle $i$ keeps decelerating and vehicle $k$ keeps accelerating, vehicle $i$ exits the constraint arc at $t_2=5.3$ $s$, when the second unconstrained arc starts.

\subsection{Case 3: linearly increasing acceleration of leading vehicle}
In this case, we consider an increasing acceleration profile of the leading vehicle $k$ with a negative initial acceleration. With the same initial speed setup, vehicle $i$ hits the constrained arc around similar time at $t_1=3.0$ $s$. However, since the speed of vehicle $k$ keeps reducing until $t=10.0$ $s$, vehicle $i$ has to decelerate for a longer time to keep the minimum safe distance with vehicle $k$. After $t=10.0$ $s$, vehicle $i$ is not able to leave the constrained arc due to the following reason: vehicle $i$ needs higher acceleration to meet the pre-defined final time, however, the acceleration of vehicle $i$ is limited by the acceleration of vehicle $k$ due to safety constraint. In case 3, we see that vehicle $i$ keeps the minimum safe following distance with vehicle $k$ until the time vehicle $k$ exits the control zone. However, if everything remains unchanged while vehicle $k$ decelerates harder, it is foreseeable that the final time of 26 $s$ is not feasible under current scenario settings.

\begin{figure} [ptb]
	\centering
	\includegraphics[width=0.5\textwidth]{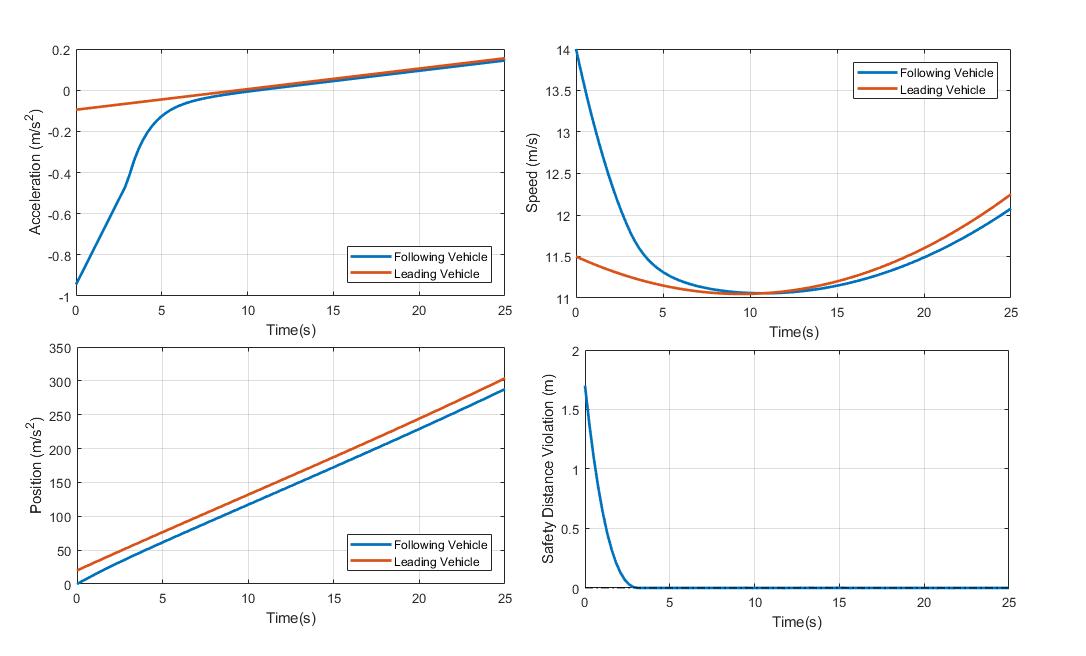}   
	\caption{Optimization results for case 3.} 
	\label{fig:case3}
\end{figure}

%%%%%%%%%%%%%%%%%%%%%%%%%%%%%%%%%%%%%%%%%%%%%%%%%%%%%%%%%%%%%%%%%%%%%%%%%%%%%%%%
%%%%%%%%%%%%%%%%%%%%%%%%%%%%%%%%%%%%
%SECTION V: CONCLUDING REMARKS AND DISCUSSION
%%%%%%%%%%%%%%%%%%%%%%%%%%%%%%%%%%%%
\section{Concluding Remarks and Discussion}
%In earlier work, a decentralized optimal control framework was established for coordinating online CAVs in merging roadways, urban intersections, speed reduction zones, and roundabouts. Hamiltonian analysis was applied to derive an analytical solution that minimizes the $L^2$-norm of the control input. However, the analytical solution did not consider the rear-end collision safety constraint. 
In this paper, we derived a closed-form analytical solution that includes the rear-end safety constraint in addition to the state and control constraints. We augmented the double integrator with an additional state corresponding to the distance of a vehicle from its preceding vehicle. Thus, we included the rear-end collision avoidance constraint as a state constraint while allowing the safe distance between the vehicles to be a function of the vehicle's speed. 
The proposed framework is limited to the lower-level individual vehicle operation control, which did not consider the upper-level vehicle coordination problem that designates the sequence that each CAV crosses the merging zone. Ongoing work considers the upper-level problem that results in maximizing the throughput of the intersection and satisfies collision avoidance constraints inside the merging zone. 
%We analyzed in the previous work that there is a trade-off between the objectives of time minimization and energy minimization, further analysis with adding time consideration in the problem formulation will be conducted. 
While the potential benefits of full penetration of CAVs to alleviate traffic congestion and reduce fuel consumption have become apparent, different penetrations of CAVs can alter significantly the efficiency of the entire system. Therefore, future research should investigate the implications of different penetration of CAVs.

%\section{ACKNOWLEDGMENTS}

%%%%%%%%%%%%%%%%%%%%%%%%%%%%%%%%%%%%%%%%%%%%%%%%%%%%%%%%%%%%%%%%%%%%%%%%%%%%%%%%
\bibliographystyle{IEEEtran}
\bibliography{TCST_references}

%\begin{IEEEbiography}[{\includegraphics[width=1.1in,height=1.25in,clip,keepaspectratio]{bio/andreas.jpg}}]{Andreas A. Malikopoulos}
% (M2006) received a Diploma in Mechanical Engineering from the National Technical University of Athens, Greece, in 2000. He received M.S. and Ph.D. degrees from the Department of Mechanical Engineering at the University of Michigan, Ann Arbor, Michigan, USA, in 2004 and 2008, respectively.
%
%He is currently the Deputy Director of the Urban Dynamics Institute and an Alvin M. Weinberg Fellow with the Energy \& Transportation Science Division at Oak Ridge National Laboratory (ORNL). Before joining ORNL he was a Senior Researcher with General Motors Global Research \& Development, conducting research in the area of stochastic optimization and control of advanced propulsion systems. His research spans several fields, including analysis, optimization, and control of complex systems; decentralized systems; and stochastic scheduling and resource allocation problems. The emphasis is on applications related to energy, transportation, and operations research.
%\end{IEEEbiography}

%\begin{IEEEbiography}
%[{\includegraphics[width=1.1in,height=1.25in,clip,keepaspectratio]{bio/andreas.jpg}}]
%{Liuhui Zhao}
%
%\end{IEEEbiography}

\end{document}